
                             \def\flech{\rightarrow}

   \NoBlackBoxes \documentstyle{amsppt}

   \magnification=1100
    \loadbold
\topmatter
 \title H\"{o}lder continuity of energy minimizer maps between riemannian
 polyhedra \endtitle

\leftheadtext{Taoufik Bouziane}

\rightheadtext{H\"{o}lder continuity of energy minimizer maps }

\redefine\headmark#1{}
                          \author Taoufik Bouziane\endauthor

     \address  Taoufik BOUZIANE, The international centre for theoretical physics mathematic section,
      strada costiera 11, 34014 Trieste Italy  \endaddress

\email tbouzian\@ictp.trieste.it, btaoufik73\@hotmail.com\endemail

\abstract

       The goal of the present paper is to establish some kind of regularity
       of an energy minimizer map between Riemannian polyhedra. More
       precisely, we will show the h\"{o}lder
      continuity of local energy minimizers between Riemannian
      polyhedra with the target spaces without focal points. With this new result,
       we also complete our existence theorem obtained in [5], and consequently
      we generalize completely, to the case
      of target polyhedra without focal points (which is weaker geometric condition than
      the nonpositivity of the curvature), the Eells-fuglede's existence and regularity
      theorem [12, chapters 10, 11] which is the new version of the famous Eells-Sampson's theorem [13].
      \endabstract


\date  \enddate



   \endtopmatter


               \document

     \head{ 0. Introduction.}\endhead

     It is well known that the problems dealing with the existence of
     energy minimizing maps are related to those of local or global
     regularity. For example, Eells and Sampson [13] proved that
     every free homotopy class of maps between smooth Riemannian
     manifolds, if the target manifolds are of nonpositive sectional
     curvature, has an energy minimizer which is smooth. Gromov
     and Schoen [15] extended the Eells-Sampson's results to
     the case when the target spaces are Riemannian polyhedra and
     obtained Lipschitz continuous energy minimizers, while
     Korevaar and Schoen [21] dropped  the polyhedral
     restriction on the target spaces permitting its to be any geodesic
     spaces. Later Eells and Fuglede [12, chapters 10, 11] proved the existence
      of h\"{o}lder continuous energy minimizers between Riemannian polyhedra with the assumption that the
     target polyhedra are of nonpositive curvature in the sense of
     Alexandrov [2]. Note that the Riemannian polyhedra are very interesting
     examples  as singular spaces, being harmonic spaces (in the sense of
     Brelot, see [12] ch 2) and providing a several examples as
     smooth Riemannian manifolds, triangulable Lipshitz manifolds,
     Riemannian orbit spaces, singular analytic spaces, stratified
     spaces, etc...

     In our turn, in [5]  we expanded the Eells-Fuglede's existence
     theorem to the case of the target polyhedron without focal
     points in the sense of [4], but the geometric arguments developed therein
     did not permit to us to tell something on the local
     regularity. In fact, we were interested in the Riemannian
     polyhedra without focal points  because this class of polyhedra
     is wider then the class of those of nonpositive Alexandrov's curvature
      even if the polyhedra are  smooth. Indeed, a geodesic space of nonpositive
     curvature is  always without focal points (cf. [4]) while, Gulliver [6] has
     shown that there are manifolds without focal points
     of both signs of sectional curvature.

      The goal of the present paper is to show the h\"{o}lder
      continuity of local energy minimizers between Riemannian
      polyhedra with the target spaces without focal points. With this
      result, we also complete our existence theorem obtained in [5], and
      consequently, we generalize completely, to the case
      of target polyhedra without focal points,
      the last version of the existence and regularity
      theorem,  due to Eells and fuglede [12, chapters 10, 11].
      Remark that, if both the source and the target of the maps are smooth,
      the theorem is due to Xin (cf. [29]).
      The methods we will use for solving this problem
      are (far) different from those used in the smooth case because the aim is also to cover
      the singular case. Consequently, we will follow the Eells-Fuglede's
      spirit, but  several difficulties arise because of our  weaker geometric
      condition (the absence of the focal points).

      The paper is organized as the following. In Section 1, we use
      the absence of the focal points in a Riemannian polyhedron to
      produce a strong convexity property of the square of the
      distance function. We note that establishing this geometric property was
      quite difficult compare to the Eells-Fuglede's case where
      the strong convexity of the square of the
      distance function is a direct consequence of the nonpositivity
      of the curvature. In section 2, we establish the h\"{o}lder continuity of energy minimizer maps
      between Riemannian polyhedra. Section 3 is devoted to
      the application of the established regularity to our
      existence theorem obtained in [5]. To close the paper, and for the sake of completeness,
       an annex  containing an overview of recent metric geometry, Riemannian
       polyhedra, Energy of map etc..., has been included with references
        for all the results stated.

        \subhead{Acknowledgements}\endsubhead

 I would like to thank Professor J.
Eells for bringing to my attention this problem, and for helpful
discussions. I owe special thanks to Professor A. Verjovsky for
his encouragement and constant support over the years.

\head { 1.  Riemannian polyhedra without focal points. }\endhead

This section is devoted to the study of the convexity of the
square of the distance function. More precisely we will firstly
investigate the case of a simply connected smooth Riemannian
manifold without focal points. Secondly, we will use the result
obtained in the smooth case to show that the square of the
distance function in simply connected Riemannian polyhedra without
focal points in sense of [4], is {\it strongly convex} (see the
definition below).

\subhead { 1.1 Smooth Riemannian manifolds without focal
points.}\endsubhead

The aim of this paragraph is to use the absence of the focal
points in a Riemannian manifold to show that the square of the
distance function from any fixed point is in some sense strongly
convex. Usually, in the literature we talk about the convexity of
the square of the distance function and never about its strong
convexity. For that reason and for our interests we decide to show
this property here.

Let $M$ denote a simply connected complete smooth Riemannian
manifold, $t\mapsto \sigma(t)$ a geodesic and $p\in M$ a point not
belonging to $\sigma$. We have the following proposition.

\proclaim{ Proposition 1.1}

If the manifold $M$ is compact without focal points, then the
square of the distance function from the point $p$ is strongly
convex, that is,

for every geodesic $\sigma : [0,1] \flech  M$, there exists a
positive constant $a$ such that:

$$d^2(p,\sigma(s))\le
(1-s)d^2(p,\sigma(0))+ s d^2(p,\sigma(1))- a s
(1-s)d^2(\sigma(0),\sigma(1)).$$

\endproclaim

For the proof of the proposition we need the following lemma.

\proclaim{ Lemma 1.2}

Let $M$ , $\sigma$ and $p$ be as in Proposition 1.1. Then there
exists a positive constant $c$ such that:
$$\frac{d^2}{ds^2}d^2(p,\sigma(s))\ge c$$

\endproclaim

\demo{Proof of Lemma 1.2}

Let us considering a geodesic variation joining the point $p$ and
the geodesic $\sigma$ as follows:

Let $\gamma : (s,t)\in \Bbb R\times [0,1]\mapsto \gamma(s,t)$ be
an $M$-valued map such that, $\gamma(.,1)=\sigma$ and for every
fixed $s$, $\gamma(s,.)$ is minimal geodesic connecting $p$ to
$\sigma(s)$.

Set $T:= \gamma_*\frac{\partial}{\partial t}$ (the direct image by
$\gamma$ of the vector field $\frac{\partial}{\partial t}$) and
$V:= \gamma_*\frac{\partial}{\partial s}$. Note that the vector
field $V$ is {\it Jacobi field} [24] with $V(0)=0$. By a direct
computation we obtain:

$$\frac{d^2}{ds^2}d^2(p,\sigma(s))=\int ^1_0 [\langle \nabla_TV,\nabla_TV\rangle -
\langle R(T,V)T,V \rangle] dt,$$ where $\nabla$ denote the {\it
symmetric Riemannian connection} of $M$ (relative to its given
Riemannian metric) and $R$ is its associate {\it curvature tensor}
[10].

 Recall that any such Jacobi field $V$ can be decomposed as follows:
$$V= V_\bot+aT+btT,$$ where $V_\bot$ is a Jacobi field
perpendicular to the geodesic $\gamma(s,.)$  and $T$ is the unit
tangent vector along $\gamma(s,.)$. Thus we obtain:

$$\frac{d^2}{ds^2}d^2(p,\sigma(s))=\int ^1_0 [\langle\nabla_TV_\bot ,\nabla_TV_\bot
\rangle- \langle R(T,V_\bot)T,V_\bot \rangle ] dt +b^2.$$

The manifold $M$ is assumed without focal points, consequently, on
one hand  the right term $\int ^1_0 [ \langle \nabla_TV_\bot
,\nabla_TV_\bot
 \rangle -\langle R(T,V_\bot)T,V_\bot\rangle] dt$ is strictly positive, on the other hand
 the real number $b$ is nonzero (because $Y=bt T$ is a nontrivial Jacobi field vanishing at $0$ $|Y(t)|$
  is strictly increasing, cf. [26], and $M$ is compact); thus the lemma
is thereby proved. \hfill $\square$

\enddemo

Next, we are going to give the proof of Proposition 1.1.

\demo{Proof of Proposition 1.1}

It is well known, in the theory of convex functions, that if we
are looking for some convexity property on a simply connected
Riemannian smooth manifold of a given function it is enough to
show the midconvexity of the relevant property. So, following this
tradition, we will show:

for given parametrization of the geodesic $\sigma$ we have for
(closed to $0$) $t\ge 0$:
$$d^2(p,\sigma(0))\le
\frac{1}{2}d^2(p,\sigma(-t))+ \frac{1}{2} d^2(p,\sigma(t))- a
t^2.$$

From Lemma 1.2, there is a positive constant $c$ such that
$\frac{d^2}{dt^2}d^2(p,\sigma(t))_{|t=0}\ge c$. Then, we obtain
for all $\epsilon>0$ the following inequality:
$$\frac{\frac{d^2(p,\sigma(t))-d^2(p,\sigma(0))}{t} -
\frac{d^2(p,\sigma(0))-d^2(p,\sigma(-t))}{t}}{t}
  \ge c - \epsilon,$$ which ends the proof of the proposition.\hfill $\square$

\enddemo

\proclaim{Remarks 1.3}

Proposition 1.1 is also valid if the manifold $M$ is the universal
cover of a compact Riemannian manifold without focal points.

\endproclaim

\subhead{1.2 Complete Riemannian polyhedra without focal
points}\endsubhead

The goal of this paragraph is to show that in the universal cover
(so simply connected) of complete compact Riemannian polyhedron
without focal points, the square of the distance function from any
fixed point is strongly convex. Recall that in [5], it is already
shown that in a simply connected locally compact Riemannian
polyhedron the square of the distance function is (just) convex.
For simplicity of statements we shall require that, our Riemannian
polyhedra are simplexwise smooth. But the results of this
paragraph are also valid with mostly the same proofs if the
Riemannian polyhedra are just Lip. \vskip .2 cm

Let $(X,d_X,g)$ be a Riemannian polyhedron endowed with
simplexwise Riemannian metric $g$ and $(K,\theta )$ a fixed
triangulation (cf. the annex 3.2).

Recall that for each point $p\in X$ (or $\theta(p)\in K$), there
are well defined notions, the tangent cone over $p$ denoted $T_pX$
and the link over $p$ noted $S_pX$ which generalizes respectively
the tangent space and the unit tangent space  if $X$ is also
smooth manifold (see the annex 3.1).

Now, we state the main result of this section.

\proclaim{Theorem 1.4}

Assume that the  polyhedron $X$ is compact and without focal
points. Let $\tilde X$ denote its universal cover and $p$ a point
of $\tilde X$. Then for every geodesic $\sigma : I\subseteq \Bbb R
\flech \tilde X$, the square of the distance function from the
point $p$ to the geodesic $\sigma$ is strongly convex, that is,

 there exists a positive constant $c$ depending only on the
polyhedron $X$ such that:

$$d^2(p,\sigma(s))\le
(1-s)d^2(p,\sigma(0))+ s d^2(p,\sigma(1))- c s
(1-s)d^2(\sigma(0),\sigma(1)).$$

\endproclaim

Before giving the proof of Theorem 1.4, let us start by the
following remarks:

\proclaim{Remarks 1.5}

\roster

\item
 Alexander and Bishop [1] have shown that a simply connected
complete locally convex geodesic space is globally convex. Thus,
following the same argument than the one used by Alexander-Bishop,
to prove Theorem 1.4, it is only required to show that every point
$x\in \tilde X$ admits an open convex neighborhood $U_x$. In other
terms, we just need to show the following: For every $x\in \tilde
X$ there is an open neighborhood $U_x$ such that, every geodesic
$\sigma$ with end points in $U_x$ belongs to $U_x$ and the
function $L: t\mapsto d^2(x,\sigma(t))$ is strongly convex.

\vskip .3 cm

\item It is shown in [5] that if $Y$ is complete simply connected
Riemannian polyhedron without focal points then, for every
geodesic $\sigma\subset Y$, the function $L: t\mapsto
d^2(x,\sigma(t))$ is continuous and convex.

\endroster

\endproclaim

\demo{Proof of Theorem 1.4}

Let $(X,g,d)$ be a compact polyhedron, $(\tilde X,\tilde g, \tilde
d)$ its universal cover (so it is complete and simply connected)
and $(K,\theta)$ a triangulation of $\tilde X$. The fundamental
group $\pi_1(X)$ acts isometrically and simplicially on $\tilde X$
thus there exists  a compact set $\tilde F\subset \tilde X$ called
a {\it fundamental domain} of $\pi_1(X)$ whose boundary $\partial
\tilde F$ has measure $0$ and each point of $\tilde X$ is
$\pi_1(X)$-equivalent either to exactly one point of the interior
of $\tilde F$ or to at least one point of $\partial \tilde F$. The
fact that $X$ is compact implies that the compact $\tilde F$ can
be obtained as a suitable finite union of maximal simplexes of
$\tilde X$.

In the following we will omit the homeomorphism of the (some)
triangulation in our notations and so we will not do any
distinction between the simplexes of $\tilde X$ and the simplexes
of $K$. By the first remark of 1.5 it is only required to show the
locale strong convexity (in our sense) and for points belonging to
the fundamental domain. So there are two cases to investigate, the
first one is when the point $p$ is in the topological interior of
some maximal simplex of the fundamental domain $\tilde F$ and the
second one is when the point $p\in \tilde F$ is vertex (to the
triangulation $(K,\theta)$).

Suppose that $p$ is in the interior of the maximal simplex
$\Delta$. Then there exists a positif reel $r_p>0$ such that the
open ball $B(p,r_p)$ with center $p$ and ray $r_p$ is contained in
$\Delta$. Thanks to the Riemannian metric $g_{\Delta}$, the open
Ball $B(p,r_p)$ can be thought of as sub-manifold of some compact
simply connected smooth Riemannian manifold endowed with the
Riemannian metric $g_{\Delta}$. Take now a geodesic $\sigma$ with
end points in $B(p,r_p)$ then by the second remark of 1.5 it is
contained in the ball $B(p,r_p)$.  The polyhedron $\tilde X$ is
without focal points so the neighborhood (sub-manifold) $B(p,r_p)$
is without focal points too. Thus, by Proposition 1.1, the
function $L$ is strongly convex for every geodesic $\sigma$
contained in $B(p,r_p)$.

Now, look at the case when $p$ is vertex of $\tilde F$. Let $r_p$
be a positif reel such that the open ball $B(p,r_p)$ is included
in the open star $st(p)$ of $p$ (see, Annex 3.1). Let $\sigma :
[a,b] \flech \tilde F$ be a geodesic of $B(p,r_p)$ and let
$\bigcup_i \Delta_i^o$ (finite union because is locally compact)
denote the star of $p$. We know that there is a subdivision
$t_0=a,t_1,...,t_n=b$ such that each restriction
$\sigma|_{[t_i,t_{i+1}]}\subset \Delta_i$ is a geodesic in  sense
of smooth Riemannian geometry. So by Proposition 1.1 and the
second remark of 1.5, the question about the strong convexity of
the function $L(t)=d^2(p,\sigma(t))$ is asked when $\sigma$
transits from a simplex $\Delta_i$ to a simplex $\Delta_{i+1}$
i.e. at the points $t_i$.

Reparameterizing the geodesic $\sigma$ and suppose that for fixed
$j$, $t_j = 0$ and that for small $\epsilon>0$ the geodesic
segment $\sigma|_{[-\epsilon,0]}$ is included in the maximal
simplex $\Delta_1$ and $\sigma|_{[0,\epsilon]}$ is included in the
maximal simplex $\Delta_2$ (by taking $\epsilon$ small enough).

By Remark 1.5 the function $L(t)=d^2(p,\sigma(t))$ is continuous
then, as it is used in the theory of convex functions, it is only
required to show the strong midconvexity of the function $L$ i.e.
$$\text{ for every $t\in [0,\epsilon]$, $L(0)\le\frac{1}{2}L(-t)+
\frac{1}{2}L(t)-ct^2.$}$$ Every maximal simplex $\Delta_i$ from
the star of $p$ is thought of as a cell of some smooth compact
simply connected Riemannian manifold $(M_i,g_i)$ without focal
points. Recall that there is an exponential function
(diffeomorphism) defined from the tangent bundle $TM_i$ of each
manifold $M_i$ to $M_i$.
 Let us now consider the two exponential maps
 $\exp_1:T_pM_1\rightarrow M_1$ and $\exp_2:T_pM_2\rightarrow M_2$.
 The spaces $T_iM$, $i=1,2$ can be assimilated to an euclidean space
 of dimension lower or equal the dimension of the polyhedron $X$. Consequently,
 there is an isometry $I: T_pM_1\rightarrow T_pM_2$. Then, thanks to
 the diffeomorphism $\exp_2 \circ I \circ \exp_1^{-1}$ (and its inverse), we can
 compare (point by point) the associate distance functions $d_1, d_2$ to the
 Riemannian metrics $g_1, g_2$ in the following sense:
 $$\text{ for every point $q\in M_1$, compare $d_1(p,q)$ and
 $d_ 2(p,\exp_2 \circ I \circ \exp_1^{-1}(q))$.}$$

With this possibility, comparing the distance functions $d_1$ and
$d_2$, we can suppose for example that $d_1\le d_2$ for some $t\in
[-\epsilon,\epsilon]$.

By the definition of the distance function $\tilde d$ of the
polyhedron $\tilde X$, we have: $\tilde d(p,\sigma(0)\le
d_1(p,\sigma(0))$. Now consider the concatenation $\sigma_1\subset
M_1$ of the two geodesics $\sigma_{|[-t,0]}$ and $\exp_1\circ
I^{-1}\circ \exp_2^{-1}(\sigma_{[0,t]})$. Then by the fact that
the polyhedron $\tilde X$ is without focal points and the distance
comparison hypothesis, $\sigma_1$ is minimal geodesic and the
function $t\mapsto d_1^2(p,\sigma_1)$ is  strongly convex i.e.

$${\tilde d}^2(p,\sigma(0))\le d_1^2(p,\sigma(0))\le \frac{1}{2} d_1^2(p,\sigma(-t))
+\frac{1}{2}d_1^2(p,\exp_1\circ I^{-1}\circ
\exp_2^{-1}(\sigma(t))) - c_1 t^2.$$

But we have supposed that $d_1\le d_2$ at $t$ so
$d_1^2(p,\exp_1\circ I^{-1}\circ \exp_2^{-1}(\sigma(t)))\le
d_2(p,\sigma(t))$ which leads to:
$$\tilde d^2(p,\sigma(0))\le d_1^2(p,\sigma(0))\le \frac{1}{2} d_1^2(p,\sigma(-t))
+\frac{1}{2}d_2^2(p,\sigma(t)) - c_1 t^2.$$ Now if we take $c=\inf
(c_1,c_2)$ we obtain:
$$\text{$\tilde d^2(p,\sigma(0))\le \frac{1}{2} d_1^2(p,\sigma(-t))
+\frac{1}{2}d_2^2(p,\sigma(t)) - c t^2.$ ($\ast$)}$$ But in the
interior of each $\Delta_i$, we have $d_i=\tilde d$ (may be we
should to take an $r_p$ smaller) so we have: $$\text{ $\tilde
d^2(p,\sigma(0))\le \frac{1}{2} \tilde d^2(p,\sigma(-t))
+\frac{1}{2}\tilde d^2(p,\sigma(t)) - c t^2$   .}$$

For ending the proof, just remark firstly that for every $t\in
[-\epsilon,\epsilon]$, the two distance functions $d_1,d_2$ are
comparable (in the above sense). Secondly the inequality $(\ast)$
is symmetric in $d_1$ and $d_2$ (because if $d_2\le d_1$ just
inverse the role of $d_1$ and $d_2$ in the proof and we obtain the
same inequality). Thirdly we can choose the $c$ uniformly because
$\tilde X$ is locally compact and $\tilde F$ is compact. Finally
the polyhedron $\tilde X$ is simply connected so by Remarks 1.5
(the first remark) the local strong convexity of the square of the
distance function becomes global for a constant $c$ depending only
on the fundamental domain $\tilde F$ and consequently it depends
only on the polyhedron $X$.

\hfill $\square$

\enddemo

Before ending this section, we mention the following remarks:

\proclaim{Remarks 1.6}

\roster

\item Theorem 1.4 is also valid if the space $\tilde X$ is compact
simply connected Riemannian polyhedron without focal points (not
necessary the universal cover of a compact Riemannian polyhedron),
or it is a Riemannian polyhedron with bounded geometry in sense of
[6].

\vskip .3cm

\item When the polyhedron $X$ is complete of nonpositive curvature
(in the sense of Alexandrov) and not necessary compact, the
constant $c$ in Theorem 1.4 is equal to $1$.
\endroster

\endproclaim

\head{2. H\"{o}lder continuity}\endhead

In the present section we will discuss a kind of regularity of an
energy minimizer map between Riemannian polyhedra. More precisely
we will take a two Riemannian polyhedra $(X,g)$ and $(Y,h)$ of
dimensions $m$ and $n$, with $X$ admissible and simplexwise
smooth, and $Y$ compact and without focal points ; and we ask the
question: what level of regularity of a given locally energy
minimizing map $\varphi :X\rightarrow \tilde Y$ ($\tilde Y$ is the
universal cover of $Y$) can we have? The best answer we obtain is
the following theorem.

\proclaim{ Theorem 2.1}

Let $X$ and $Y$ be Riemannian polyhedra. Suppose that $X$ is
admissible and simplexwise smooth and  $Y$ is compact without
focal points. Then, every locally energy minimizing map $\varphi
:X\rightarrow \tilde Y$, where $\tilde Y$ is the universal cover
of $Y$, is H\"{o}lder continuous.

\endproclaim

Recall that in our context, a map $\varphi :X\rightarrow \tilde Y$
is H\"{o}lder continuous if there is a H\"{o}lder continuous map
which is equals to $\varphi$ almost everywhere in $X$.

To prove Theorem 2.1, we will adapt to our frame the arguments
used by Eells and Fuglede in [12] where they proved a similar
theorem but in the case of the target polyhedron a  complete
Riemannian polyhedron of nonpositive curvature. However, the
difficulties in our case come firstly, from the fact that we are
considering the Riemannian polyhedra as geometric habitat and
where, in general, we can not use the second differential
calculus. Secondly, difficulties also arise from the fact that the
strong convexity of the square of the distance function
established in Section 1 is quite weaker than the convexity
property  used by Eells-Fuglede  in [12] and it is also optimal in
our case.

Now let us begin with some lemmas which will be needed for the
proof of Theorem 2.1.

\proclaim{Lemma 2.2}

Let $Y$ be a compact Riemannian polyhedron without focal points.
Let $\varphi \in W^{1,2}_{loc}(X,\tilde Y)$ be locally energy
minimizing map, where $(\tilde Y,d_{\tilde Y})$ is the universal
cover of $Y$. Then $\varphi$ is essentially locally bounded (i.e.
$\varphi=\tilde \varphi, a.e.$ and $\tilde \varphi$ is bounded)
and we have for any $q\in \tilde Y$:

\roster

\item The functions $v_q: x\mapsto d_{\tilde Y}(q,\varphi(x))$ and
$v_q^2:x\mapsto d_{\tilde Y}^2(q,\varphi(x))$ from
$W^{1,2}_{loc}(X)$ are weakly subharmonic.

\item $E(v_q^2,\lambda)\le -2 \frac{1}{c} \int_X e(\varphi)\lambda
d\mu_g$ for every $\lambda\in W^{1,2}_c(X)\bigcap L^\infty(X)$,
$\lambda\ge 0$ which we can write in weak sense: $\triangle v_q^2
\ge  2 \frac{1}{c} e(\varphi)$, with $c$ is the constant of the
  strong  convexity of the square of the distance function, where
$E(v_q^2,\lambda):=\int_X c_m \langle\nabla\lambda, \nabla
v_q^2\rangle d \mu_g$ with $c_m = \frac{\omega_m}{m+1}$ and,
$\omega_m$ being the volume of the unit ball in $\Bbb R^m$ and
$\nabla$ and $\langle.,.\rangle$ denote respectively the gradient
operator and the inner product, defined a.e in $X$ .
\endroster

\endproclaim

\demo{Proof Lemma 2.2}

Following the same idea used by Eells-Fugled (Lemma 10.2 in [12]),
an idea used by Jost [19] and before used in [20], we will compare
the map $\varphi$ as in the Lemma 2.2 with maps obtained by
pulling $\varphi(x)$ towards a given point $q\in \tilde Y$.

The map $\varphi: X\rightarrow \tilde Y$ from the space
$W^{1,2}_{loc}(X,\tilde Y)$ is locally energy minimizing, so $X$
can be covered by relatively compact domains $U\subset X$ for
which $E(\varphi_{|U})\le E(\psi_{|U})$ for every map $\psi\in
W^{1,2}_{loc}(X,\tilde Y)$ such that $\psi=\varphi$ a.e. in
$X\setminus U$.

Note $v=v_q: x\mapsto d_{\tilde Y}(q,\varphi(x))$,
$v^2=v_q^2:x\mapsto d_{\tilde Y}^2(q,\varphi(x))$ and set $d$
refereing to the distance function $d_{\tilde Y}$ in $\tilde Y$.
All the properties we want to show are local so we can suppose
that $X$ is compact and it does not alter the results of the
lemma.

By Theorem 1.4 the square of the distance function is strongly
convex so there exists a constant $c>0$ (depending only on the
polyhedron $Y$) such that for every geodesic $\sigma$ arc length
parameterized we have:

$$d^2(p,\sigma(s))\le
(1-s)d^2(p,\sigma(0))+ s d^2(p,\sigma(1))- c s
(1-s)d^2(\sigma(0),\sigma(1)).$$

There are two case to investigate, the first one is when $c\ge 1$
and the second one is when $c<1$.

In the first case $c\ge1$, the strong convexity of the square
distance function implies the following:

$$d^2(p,\sigma(s))\le
(1-s)d^2(p,\sigma(0))+ s d^2(p,\sigma(1))-  s
(1-s)d^2(\sigma(0),\sigma(1)),$$ which join the case of
nonpositive curvature and so by Eells-fuglede's results ([12] ch
10) the Lemma 2.2 follows.

Now look at the second case when $0<c<1$. Let us first consider
the case when $\lambda$ is lipschitz map on $X$, $0\le \lambda\le
1$ and the support of $\lambda$ noted $supp \lambda$ is subset of
some domain $U$ as we mentioned in the beginning of the proof. Let
$\gamma_x$ denote the minimal geodesic in the space $\tilde Y$
joining the point $\gamma_x(0)=\varphi(x)$ to $\gamma_x(1)=q$ for
$x\in X$. Define a map (the pulling of $\varphi$) $\varphi_\lambda
:X\rightarrow \tilde Y$ by
$\varphi_\lambda(x)=\gamma_x(\lambda(x))$, for $x\in X$.

The function $\varphi_\lambda$ is $L^2(X, \tilde Y)$ because the
geodesic $\gamma_x$ varies continuously  with its end point
$\varphi(x)$ (The space $\tilde X$ is without conjugate points,
see [5]) and $d(\varphi_\lambda(x),q)\le d(\varphi(x),q)$ with
$d(\varphi(.),q)\in L^2(X)$ (since $\varphi\in L^2(X,\tilde Y)$).

The strong convexity of the square of the distance function gives,
for $x,x'\in X$,

$$d^2(\varphi(x),\varphi_\lambda(x'))\le
   (1-\lambda(x'))d^2(\varphi(x),\varphi(x'))+\lambda(x')d^2(\varphi(x),q)-$$
   $$c \lambda(x')(1-\lambda(x'))d^2(\varphi(x'),q)$$ and
\vfill
 \pagebreak

   $$d^2(\varphi_\lambda(x),\varphi_\lambda(x'))\le
   (1-\lambda(x))d^2(\varphi(x),\varphi_\lambda(x'))+\lambda(x)d^2(\varphi_\lambda(x'),q)-$$
   $$c \lambda(x)(1-\lambda(x))d^2(\varphi(x),q).$$
Combining the two inequalities and inserting
$d(\varphi_\lambda(x'),q)=(1-\lambda(x'))d(\varphi(x'),q)$, we
obtain,
$$d_\lambda^2-d^2\le
-[\lambda(x)+\lambda(x')-\lambda(x)\lambda(x')]d^2 +
(1-\lambda(x))\lambda(x')v^2(x)- c\lambda(x)(1-\lambda(x))v^2(x)$$
$$- c(1-\lambda(x))\lambda(x')(1-\lambda(x'))v^2(x')+ \lambda(x)(1-\lambda(x'))^2 v^2(x'),$$
with $d_\lambda=d(\varphi_\lambda(x),\varphi_\lambda(x'))$ and
$d=d(\varphi(x),\varphi_(x'))$.

Now, we will pull closely the map $\varphi$ towards the point $q$
such that $\lambda\ge \frac{1}{1+c}$. Under this assumption
($\lambda\ge \frac{1}{1+c}$), we have for every $x,x'\in X$,
$\lambda(x)(1-\lambda(x'))\le c\lambda(x)(1-c\lambda(x'))$ and
$\lambda(x')(1-\lambda(x))\le c\lambda(x')(1-c\lambda(x))$. Taking
in account these inequalities we obtain,

$$d_\lambda^2-d^2\le
-[\lambda(x)+\lambda(x')-\lambda(x)\lambda(x')]d^2 - c
(\lambda(x)-\lambda(x')) (v^2(x)-v^2(x'))$$
$$\text{$+$ $|O(\lambda^2(x))|$ $+$ $|O(\lambda(x')\lambda(x))|$ $+$
 $|O(\lambda(x)\lambda^2(x'))|$ .}$$

For given $\lambda$ a compact supporting positive lipschitz
function, replace $\lambda$ with $t\lambda$, for any
$0<t<\epsilon$, this leads to,

$$d_{t\lambda}^2-d^2\le
-[t\lambda(x)+t\lambda(x')-t^2\lambda(x)\lambda(x')]d^2 - c
t(\lambda(x)-\lambda(x')) (v^2(x)-v^2(x'))$$
$$\text{$+$ $\epsilon^2|O(\lambda^2(x))|$ $+$ $\epsilon^2|O(\lambda(x')\lambda(x))|$ $+$
 $\epsilon^3|O(\lambda(x)\lambda^2(x'))|$ \hskip .5cm $(\ast)$.}$$

 Observe that
 $$\text{$ E({\varphi_{t\lambda}}_{|U})-E(\varphi_{|U})=$}
 \lim_{f\in C_c(U,[0,1])} \limsup_{\epsilon\rightarrow 0}\int_U
(e_\epsilon (\varphi_{t\lambda})- e_\epsilon(\varphi) f d\mu_g,$$
where for $\epsilon>0$,
$e_\epsilon(\psi)(x)=\int_{B_X(x,\epsilon)}\frac{d_{\tilde
Y}^2(\psi(x),\psi(x'))}
 {\epsilon^{m+2}}d\mu_g(x')$, with $\psi\in L^2_{loc}(X,\tilde
 Y)$.

  Note that for $x'\in B_X(x,\epsilon)$ we have
 $|\lambda(x)-\lambda(x')|\le const. \epsilon$ and $\varphi$ is supposed locally
 energy minimizing, so by Eells-Fuglede's results
[12, Definition 9.1, Theorem 9.1 and Corollary 9.2] we deduce from
the inequality $(\ast)$,
 $$0\le E({\varphi_{t\lambda}}_{|U})-E(\varphi_{|U})
 \le -\int_U (2t\lambda- t^2 \lambda^2) e(\varphi)
d\mu_g -c c_m\int_U t \langle \nabla \lambda,\nabla v^2 \rangle
d\mu_g +$$
 $$\text{ $\int_U$ $c_m
(O(||\lambda||_{L^\infty}^2)$ $+$ $O(||\lambda||_{L^\infty}^3))
d\mu_g$ ,}$$ where $c_m$ is constant depending on the polyhedron
$X$ and, $\nabla$ and $\langle.,.\rangle$ denote respectively the
gradient operator and the inner product, defined a.e in $X$ (cf.
[21] and [12, ch 5]). Thus the function $\varphi_{t\lambda}$ is in
the space $ W^{1,2}_{loc}(X)$.

Now again replace in the last inequalities $\lambda$ with
$t^2\lambda$ divided by $t^3$ and let $t\rightarrow 0$ we obtain
for every $\lambda\in Lip^+_c(U)$,
$$0\le -\int_U
2\lambda e(\varphi) d\mu_g -c c_m\int_U \langle \nabla
\lambda,\nabla v^2\rangle d\mu_g. $$ So we infer that
 for every $\lambda\in
Lip^+_c(U)$,
$$0\le \int_U 2\lambda e(\varphi) d\mu_g \le -c c_m\int_U \langle \nabla \lambda,\nabla
v^2\rangle d\mu_g. $$ These inequalities extend to a functions
$\lambda\in W^{1,2}_c(U)\cap L^\infty(U)$, with $\lambda\ge 0$,
because any such $\lambda$ can be approximated in $W_c^{1,2}(U)$
by uniformly bounded functions in $L_c^+(U)$. Thus for $\lambda\in
W^{1,2}_c(X)\cap L^\infty(X)$, $\lambda\ge 0$ we have (on $X$),
$$0\le \int_X 2\lambda e(\varphi) d\mu_g \le -c c_m\int_X \langle \nabla \lambda,\nabla
v^2\rangle d\mu_g. $$ So we have shown the second statement of
Lemma 2.2. For the first part of the lemma, just remark that by
the last inequalities we have, for every $\lambda\in
W^{1,2}_c(X)\cap L^\infty(X)$, $\lambda\ge 0 $,
$$\int_X \langle \nabla \lambda,\nabla
v^2\rangle d\mu_g\le 0 , $$ which means that the function $v^2=
d^2(\varphi(.),q)$ is weakly subharmonic in $X$, and in particular
essentially locally bounded.

For the function  $v= d(\varphi(.),q)$, by the usual polarization
[12, page 21 (2.1)] we have for every $\lambda\in W^{1,2}_c(X)\cap
L^\infty(X)$, $\lambda\ge 0 $,
$$E(v^2,\lambda)= 2 E(v,\lambda v)- 2 c_m \int_X \lambda \langle
\nabla v, \nabla v\rangle d\mu_g.$$ Remember that by the triangle
inequality, $|v(x)-v(x')|^2\le d^2(\varphi(x),\varphi(x'))$, and
so [12, corollary 9.2] $c_m|\nabla v|^2\le e(\varphi)$. Inserting
$E(v^2,\lambda)\le -\frac{1}{c}\int_X 2\lambda e(\varphi) d\mu_g
\le 0$ in the last equality, it therefor follows, for every
$\lambda\in W^{1,2}_c(X)\cap L^\infty(X)$, $\lambda\ge 0 $,
$$ E(v,\lambda v)= \int_X c_m\langle
\nabla v, \nabla (\lambda v)\rangle d\mu_g \le
(1-\frac{1}{c})\int_X \lambda e(\varphi) d\mu_g.$$ But the
constant $c$ is supposed $<1$ so we deduce that,
$$\text{ for every $\lambda\in W^{1,2}_c(X)\cap L^\infty(X)$,
$\lambda\ge 0 $, $\int_X \langle \nabla v, \nabla (\lambda
v)\rangle d\mu_g \le 0$}.$$ Now, using the Eells-Fuglede's
arguments [12, page 183], we deduce that the function $v$ is
weakly subharmonic in X and it is essentially locally bounded too.

\hfill $\square$

\proclaim{Corollary 2.3}

Under the hypotheses of Lemma 2.2, if the map $\varphi$ is of
(global) finite energy then, $E(v^2,\lambda)\le -2 \frac{1}{c}
\int_X e(\varphi) \lambda d\mu_g$ for every $\lambda\in
W^{1,2}_0(X)\cap L^\infty(X)$, $\lambda\ge 0 $.

\endproclaim

\demo{Proof of Corollary 2.3}

By Lemma 2.2, we have, for every $\lambda\in W^{1,2}_c(X)\cap
L^\infty(X)$, $\lambda\ge 0 $, $E(v^2,\lambda)\le -2 \frac{1}{c}
\int_X e(\varphi) \lambda d\mu_g$.

By truncation, any positive function $\lambda\in W^{1,2}_0 (X)\cap
L^\infty(X)$ can be approximated in $W^{1,2}_0(X)$ by a uniformly
bounded sequence of function $\lambda_n\in W^{1,2}_c(X)\cap
L^\infty(X)$.

Now, $E(\varphi)$ is supposed $<\infty$,  so by the dominated
convergence theorem, we have $\int \lambda_n e(\varphi) d\mu_g
\rightarrow \int \lambda e(\varphi) d\mu_g$ (in fact there is a
subsequence of $(\lambda_n)$ which converges to $\lambda$
pointwise $a.e.$)

\hfill $\square$

\enddemo

\proclaim{Corollary 2.4}

Under the hypotheses of Lemma 2.2, every locally energy minimizing
map $\varphi :X \rightarrow \tilde Y$ is locally essentially
bounded.

\endproclaim

\demo {Proof of Corollary 2.4}

The Eells-fuglede's proof [12, Corollary 10.1] of the same
statement in the case where the target space $Y$ is of nonpositive
curvature, remains valid in our case.

\hfill $\square$

\enddemo

The following lemma will be necessary in the proof of the next
one, but it is also of special self interest.

\proclaim{Lemma 2.5}

Let $(\frak N, \nu)$ be a probability measure space, let $(\tilde
Y, d)$ be the universal covering of compact Riemannian polyhedron
$Y$ without focal points, and $f\in L^2(\frak N, \tilde Y)$. Then
there exists a unique center of mass $\bar{f_\nu}$, defined as the
point in $\tilde Y$ which minimizes the integral $\int_{\frak N}
d^2(f(x),q) d\nu(x)$.
\endproclaim

\demo{Proof of Lemma 2.5}

The space $\tilde Y$ is supposed without focal points,
consequently the square of the distance function is strong convex.
So, if $y_1, y_2$ are two points in $\tilde Y$ and
$y_{\frac{1}{2}}$ is their midpoint (the unique point in the
unique geodesic between $y_1$ and $y_2$ which is at equal distance
to both $y_1$ and $y_2$), then we have,

$$d^2(f(x),y_{\frac{1}{2}})\le \frac{1}{2}d^2(f(x),y_1)+ \frac{1}{2}d^2(f(x),y_2)
        -\frac{1}{4} c d^2(y_1,y_2),$$ with $0<c$ a constant
        depending on the space $Y$.
        Integrating over $\frak N$ we obtain,
        $$ \frac{1}{4} c d^2(y_1,y_2)\le \frac{1}{2}\int_{\frak N}
        d^2(f(x),y_1)d\nu(x)
        + \frac{1}{2}\int_{\frak N} d^2(f(x),y_2)d\nu(x)
        -  \int_{\frak N}  d^2(f(x),y_{\frac{1}{2}}) d\nu(x),$$
        Thus any minimizing sequence $(x_i)$ is Cauchy, in particular it converges to a unique
        limit point ($\tilde Y$ is complete) which is the unique minimizer of our
        integral.

\hfill $\square$

\enddemo

\proclaim{Lemma 2.6}

Let $(X,g)$ denote a compact admissible Riemannian polyhedron and
$(\tilde Y, d)$ be the universal covering of compact Riemannian
polyhedron $Y$ without focal points. For every measurable set
$A\subset X$ with $\mu_g(A)>0$, the meanvalue $\bar{\varphi_A}\in
\tilde Y$ over $A$ of a map $\varphi \in W^{1,2}(X,Y)$, defined as
the minimizing point in $\tilde Y$ of the integral $\int_A
d^2(\varphi(x),q) d\mu_g(x)$,
 lies in the closed convex hull of the essential image
 $\varphi(A)$.
\endproclaim

Before proving the lemma, just recall that the {\it essential
image } $\varphi(A)$ is defined as the closed set of all points
$q\in \tilde Y$ such that $A\cap \varphi^{-1}(V)$ has positive
measure for any neighborhood $V$ of $q$ in $\tilde Y$, and the
{\it closed convex hull} of a set $B\subset \tilde Y$ is defined
as the intersection of all closed convex subsets of $\tilde Y$
containing $B$.

\demo{Proof of Lemma 2.6}

The existence and uniqueness of the meanvalue point (of a map
belonging to $W^{1,2}(X,Y)$) over any measurable subset of $X$ are
immediately deduced from Lemma 2.5.

Let $C$ denote any convex set containing $\varphi(A)$, let $y\in
\tilde Y\setminus C$, we claim that there is unique point $\hat{y}
\in C$ nearest to $y$. Indeed, any minimizing sequence
$(y_i)\subset C$ for the function $d^2(y,.)$ on $C$, is Cauchy in
$\tilde Y$, by the strong convexity  of the square of the distance
function, and hence has a unique limit point $\hat{y} \in C$
(because $C$ is closed).

Now, let $z$ denote any point of $C$, consider the unique geodesic
$\sigma_{z\hat{y}}$ (because $\tilde Y$ is without conjugate
points [5]) connecting $z$ to $\hat{y}$. The point $\hat{y}$ is
the unique orthogonal projection (in sense of [5]) of the point
$y$ on the geodesic $\sigma_{z\hat{y}}$, consequently the angle at
the point $\hat{y}$ (the distance in the link of $\hat{y}$)
between the geodesics $\sigma_{z\hat{y}}$ and $\sigma_{\hat{y}y}$
is $\ge \frac{\pi}{2}$, and so, we deduce that, $d(z,y)>
d(z,\hat{y})$.

Moreover, we infer that $d(\varphi(.),y)>d(\varphi(.),\hat{y})$
$a.e$ in $A$, which rules out the possibility that any point $y\in
\tilde Y \setminus C$ can be the meanvalue of $\varphi$ over $A$.

\hfill $\square$
\enddemo

Now, we have all the ingredients to prove Theorem 2.1.

 \demo{Proof of Theorem 2.1}

Replacing  Lemma 10.2 , Lemma 10.4 and Corollary 10.1 in the
Eells-Fuglede's proof of the equivalent theorem [12, Ch10 , page
189] in the case of the target polyhedron of nonpositive
curvature, with our Lemma 2.2, Lemma 2.6 and Corollary 2.4
respectively, using the weak Poincar\'e inequality [cf. 12,
Proposition 9.1] and the fact that Sublemma 10.1, Lemma 10.3 and
Corollary 10.2 of [12] remain valid in our case, then using the
Eells-Fuglede's arguments [12, pages 189-192], we easily derive
our theorem.

\hfill $\square$
\enddemo

\enddemo

\head{3. Application} \endhead

A naturel question comes to our minds after the regularity result
of Theorem 2.1: is when or where can we apply the regularity
obtained? Thus, we will close the paper with this short section
where we will give an example of such application. The application
proposed will in some sense, complete the existence result  of
energy minimizer maps, obtained in [5]. Henceforth all polyhedra
considered are supposed simplexwise smooth.

\proclaim{ Theorem 3.1}

Let $X$ and $Y$ be compact Riemannian polyhedra. Suppose that $X$
is admissible and $Y$ is without focal points.

Then every homotopy class $[u]$ of each continuous map $u$ between
the polyhedra $X$ and $Y$ has an energy minimizer relative to
$[u]$ which is H\"{o}lder continuous.

\endproclaim

\demo{ Proof of Theorem 3.1}

Let $X$ and $Y$ be two compact Riemannian polyhedra such that $X$
is admissible and $Y$ is without focal points.

Firstly, remark that the existence part of Theorem 3.1 is already
proved in [5].

Secondly, as we showed in the proof of the existence part in [5],
if $u$ denote an energy minimizer in the class $[u]$ then it can
be covered by a map $\tilde u : \tilde X\rightarrow \tilde Y$,
where $\tilde X$ and $\tilde Y$ denote respectively the universal
covers of $X$ and $Y$, and which minimizes the energy in the class
of the equivariant maps with respect to the fundamental groups
$\pi_1(X)$ and $\pi_1(Y)$ in $W^{1,2}(\tilde X,\tilde Y)$.
Moreover $E(u)=\int_{\tilde F} e(\tilde u)$, where $\tilde
F\subset \tilde X$ denote the fundamental domain of $\pi_1(X)$.
But the universal cover $\tilde Y$ satisfies the hypothesis of
Theorem 2.1 so $\tilde u$ is h\"{o}lder continuous, and therefore
is the map $u$ (the energy minimizer map relative to the class
$[u]$). This ends the proof.

\hfill $\square$

\enddemo

                 \head { Annex. }\endhead

The annex is globally devoted to an overview concerning the
geodesic spaces, Riemannian polyhedra and the  harmonic maps on
singular spaces. The last subject was developed successively by
Gromov-Schoen [15], Korevaar-Schoen [21] [22] and Eells-Fuglede
[12]. We hope that the annex will be useful for the reader.

\subhead { 1. Geodesic spaces [2] [6] [7] [8] [14] }\endsubhead

Let $X$ be a metric space with metric $d$. A curve $c:I\flech X$
is called a {\it geodesic} if there is $v\geq 0$, called the
speed, such that every $t\in I$ has neighborhood $U\subset I$ with
$d(c(t_1),c(t_2))=v|t_1-t_2|$ for all $t_1,t_2\in U$. If the above
equality holds for all $t_1,t_2\in I$, then $c$ is called {\it
minimal geodesic}.

The space $X$ is called a {\it geodesic space} if every two points
in $X$ are connected by minimal geodesic. We assume from now on
that $X$ is complete geodesic space.

A triangle $\Delta$ in $X$ is a triple $(\sigma_1, \sigma_2,
\sigma_3)$ of geodesic segments whose end points match in the
usual way. Denote by $H_k$ the simply connected complete surface
of constant Gauss curvature $k$. A {\it comparison triangle}
$\bar{\Delta}$ for a triangle $\Delta \subset X$ is a triangle in
$H_k$ with the same lengths of sides as $\Delta$. A comparison
triangle in $H_k$ exists  and is unique up to congruence if the
lengths of sides of $\Delta$ satisfy the triangle inequality and,
in the case $k>0$, if the perimeter of $\Delta$ is
$<\frac{2\pi}{\sqrt k }$. Let $\bar{\Delta}= (\bar\sigma_1,
\bar\sigma_2, \bar\sigma_3)$ be a comparison triangle for
$\Delta=(\sigma_1, \sigma_2, \sigma_3)$, then for every point
$x\in \sigma_i$, $i=1,2,3$, we denote by $\bar x$ the unique point
on $\bar \sigma_i$ which lies at the same distances to the ends as
$x$.

Let $d$ denote the distance functions in both $X$ and $H_k$. A
triangle $\Delta$ in $X$ is $CAT_k$ {\it triangle} if the sides
satisfy the triangle inequality, the perimeter of $\Delta$ is
$<\frac{2\pi}{\sqrt k }$ for $k>0$, and if $d(x,y)\le d(\bar
x,\bar y)$, for every two points $x,y\in X$.

We say that $X$ has curvature at most $k$ and write $k_X\le k$ if
every point $x\in X$ has a neighborhood $U$ such that any triangle
in $X$ with vertices in $U$ and minimizing sides is $CAT_k$. Note
that we do not define $k_X$. If $X$ is Riemannian manifold, then
$k_X\le k$ iff $k$ is an upper bound for the sectional curvature
of $X$.

A geodesic space $X$ is called geodesicaly complete iff every
geodesic can be stretched in the two direction.

We say that a geodesic space $X$ is without conjugate points if
every two points in $X$ are connected by unique geodesic.

\subhead { 2. Orthogonality and focal point}\endsubhead

For more details on the study of focal points in geodesic space,
the reader can see [4] and [5].

\definition{2.1 Orthogonality}

$(X,d)$ will denote a complete geodesic space. Let $\sigma : \Bbb
R \flech X$ denote a geodesic and $\sigma _1:[a,b]\flech X$ a
minimal geodesic with a foot in $\sigma$ (i.e. $\sigma_1(a)\in
\sigma(\Bbb R)$).

The geodesic $\sigma_1$ is {\it orthogonal} to $\sigma$ if for all
$t\in [a,b]$, the point $\sigma_1(t)$ is locally of minimal
distance from $\sigma$.

In the case when for given geodesic $\sigma$ and a non-belonging
point $p$ there exists an orthogonal geodesic $\sigma'$ to
$\sigma$ and containing $p$, we will call the intersection point
between $\sigma$ and $\sigma'$ the {\it orthogonal projection
point} of $p$ on $\sigma$.

It is shown in [4] that, on the one hand, if the geodesic $\sigma$
is minimal then  there always exists a realizing distance
orthogonal geodesic to $\sigma$ connecting every external point
$p$ (off $\sigma$)  to $\sigma$. On the other hand, If the space
$(X,d)$ is locally compact with non-null injectivity radius and
the geodesic $\sigma$ is minimal on every open interval with
length lower than the injectivity radius, then for every point $p$
off $\sigma$ and whose distance from $\sigma$ is not greater than
the half of the injectivity radius, there exists a geodesic
joining orthogonally the point $p$ and the geodesic $\sigma$.

As corollaries, if the space $(X,d)$ is simply connected $CAT_0$
space then for given geodesic $\sigma : \Bbb R \flech X$ and an
off point $p$ there always exists a realizing distance orthogonal
geodesic from $p$ to $\sigma$. When $X$ is $CAT_k$ for positive
constant $k$ then there always exists an orthogonal geodesic to
$\sigma$ from a point $p$ whose distance from $\sigma$ is not
greater than $\pi\over {2\sqrt k}$. In these last two cases the
angle between two orthogonal geodesics (in the sense of the
definition above) is always greater than or equal to $\pi\over 2$.

\enddefinition
 \definition{2.2 Focal points}

Let $(X,d)$ denote a complete geodesic space, $\sigma :\Bbb R
\flech X$ a geodesic and $p$ a point not belonging to the geodesic
$\sigma$.

The point $p$ is said a {\it focal point} of the geodesic $\sigma$
or just a focal point of the space $X$, if there exists a minimal
geodesic variation $\tilde \sigma :]-\epsilon ,\epsilon [
\times[0,l]\flech X$ such that, if we note $\tilde \sigma
(t,s)=\sigma _t(s)$,  $\sigma _0$ is minimal geodesic joining $p$
to the point $q=\sigma (0)$ and for every $t\in ]-\epsilon
,\epsilon [$, $\sigma_t$ is minimal geodesic containing
$\sigma(t)$, with the properties:

  \roster

  \item For every $t\in ]-\epsilon ,\epsilon[$,
        each a geodesic $\sigma _t$ is orthogonal to $\sigma$.

        \item
       $\lim\limits _{t\to 0}{d(p,\sigma _t(l))\over {d(q,\sigma (t))}} =0$.
                             \endroster

 \enddefinition

 This definition was introduced in [4], as a
 natural generalization of the same notion in the smooth case.
 It is shown in the same paper that the Hadamard spaces are
 without a focal point.

 It is also shown in [5], that if $X$ is simply connected geodesic
 space without focal points then it is without conjugate points.

\subhead { 3. Riemannian polyhedra}\endsubhead

\definition{3.1 Riemannian admissible complexes ([3] [6] [7] [11] [28])}

 Let $K$ be locally finite simplicial complex, endowed with a
 piecewise smooth Riemannian metric $g$; i.e. $g$ is a family of
 smooth Riemannian metrics $g_\Delta$ on simplices $\Delta$ of $K$
 such that the restriction $g_\Delta|\Delta'=g_{\Delta'}$ for any
 simplices $\Delta'$ and $\Delta$ with $\Delta'\subset \Delta$.

 Let $K$ a finite dimensional simplicial complex which is connected
locally finite. A map $f$ from $[a,b]$ to $K$ is called a broken
geodesic if there is a subdivision $a=t_0<t_1<...<t_{p+1}=b$ such
that $f([t_i,t_{i+1}])$ is contained in some cell and the
restriction of f to $[t_i,t_{i+1}]$ is a geodesic inside that
cell. Then define the length of the broken geodesic map $f$ to be:
 $$ L(f)=\sum_{i=0}^{i=p} d(f(t_i),f(t_{i+1})).$$
 The length inside each cell being measured with respect its metric.

 Then define $\tilde d(x,y)$, for every two points $x,y$ in $K$,
 to be the lower bound of the lengths of broken geodesics from $x$
 to $y$. $\tilde d$ is a pseudo-distance.

 If $K$ is connected and locally finite, then
$(K,\tilde d)$ is length space which is geodesic space if complete
(see also [6]).

A $l$-simplex in $K$ is called a {\it boundary simplex} if it is
adjacent to exactly one $l+1$ simplex. The complex $K$ is called
{\it boundaryless} if there are no boundary simplices in $K$.

The (open) {\it star} of an open simplex $\Delta^o$ (i.e. the
topological interior of $\Delta$ or the points of $\Delta$ not
belonging to any sub-face of $\Delta$, so if $\Delta$ is point
then $\Delta^o=\Delta$) of $K$ is defined as:
$$\text{$st(\Delta^o)=\bigcup \{ \Delta_i^o : \Delta_i$ is simplex
of $K$ with $\Delta_i\supset\Delta \}$ .}$$ The star $st(p)$ of
point $p$ is defined as the star of its {\it carrier}, the unique
open simplex $\Delta^o$ containing $p$. Every star is path
connected and contains the star of its points. In particular $K$
is locally path connected. The closure of any star is sub-complex.

\hskip 0 cm

 We say that the complex $K$ is {\it admissible}, if it is dimensionally homogeneous
and for every connected open subset $U$ of $K$, the open set
$U\setminus \{ U\cap \{\text{the $(k-2)$}-\text{skeleton}\} \}$ is
connected ($k$ is the dimension of $K$)(i.e. K is
$(n-1)$-chainable).

\hskip 0 cm

Let $x\in K$ a vertex of $K$ so that $x$ is in the $l$-simplex
$\Delta_{l}$. We view $\Delta_{l}$ as an affine simplex in $\Bbb
R^l$, that is $\Delta _l =\bigcap_{i=0}^l H_i$, where
$H_0,H_1,...,H_l$ are closed half spaces in general position, and
we suppose that $x$ is in the topological interior of $H_0$. The
Riemannian metric $g_{\Delta_l}$ is the restriction to $\Delta_l$
of a smooth Riemannian metric defined in an open neighborhood $V$
of $\Delta_l$ in $\Bbb R^l$. The intersection
$T_x\Delta_l=\bigcap_{i=1}^l H_i \subset T_xV$ is a cone with apex
$0\in T_xV$, and $g_{\Delta_l}(x)$ turns it into an euclidean
cone. Let $\Delta_m\subset \Delta_l$ ($m<l$) be another simplex
adjacent to $x$. Then, the face of $T_x\Delta_l$ corresponding to
$\Delta_m$ is isomorphic to $T_x\Delta_m$ and we view
$T_x\Delta_m$ as a subset of $T_x\Delta_l$.

Set $T_xK =\bigcup_{\Delta_i\ni x} T_x\Delta_i$, we call it the
{\it tangent cone} of $K$ at $x$. Let $S_x\Delta_l$ denote the
subset of all unit vectors in $T_x\Delta_l$ and set $S_x=S_xK
=\bigcup_{\Delta_i\ni x} S_x\Delta_i$. The set $S_x$ is called the
{\it link} of $x$ in $K$. If $\Delta_l$ is a simplex adjacent to
$x$, then $g_{\Delta_l}(x)$ defines a Riemannian metric on the
$(l-1)$-simplex $S_x\Delta_l$. The family $g_x$ of Riemannian
metrics $g_{\Delta_l}(x)$ turns $S_x\Delta_l$ into a simplicial
complex with a piecewise smooth Riemannian metric such that the
simplices are spherical.

We call an admissible  connected locally finite simplicial
complex, endowed with a piecewise smooth Riemannian metric, an
{\it admissible Riemannian complex}.

 \enddefinition

\definition{3.2 Riemannian polyhedron [9], [1]}

We mean by {\it polyhedron} a connected locally compact separable
Hausdorff space $X$ for wich there exists a simplicial complex $K$
and homeomorphism $\theta : K \flech X$. Any such pair $(K,\theta
)$ is called a {\it triangulation} of $X$. The complex $K$ is
necessarily countable and locally finite (cf. [27] page 120) and
the space $X$ is path connected and locally contractible. The {\it
dimension} of $X$ is by definition the dimension of $K$ and it is
independent of the triangulation.

A {\it sub-polyhedron} of a polyhedron $X$ with given
triangulation $(K,\theta )$, is polyhedron $X'\subset X$ having as
a triangulation $(K',\theta|_{K'})$ where $K'$ is a subcomplex of
$K$ (i.e. $K'$ is complex whose vertices and simplexes are some of
those of $K$).

If $X$ is polyhedron with specified triangulation $(K,\theta)$, we
shall speak of vertices, simplexes, $i-$skeletons or stars of $X$
respectively of a space of links or tangent cones of $X$ as the
image under $\theta$ of vertices, simplexes, $i-$skeletons or
stars of $K$ respectively the image of space of links or tangent
cones of $K$. Thus our simplexes become compact subsets of $X$ and
the $i-$skeletons and stars become sub-polyhedrons of $X$.

If for given triangulation $(K,\theta)$ of the polyhedron $X$, the
homeomorphism $\theta$ is locally bi-lipschitz then $X$ is said
{\it Lip polyhedron} and $\theta$ {\it Lip homeomorphism}.

A {\it null set} in a Lip polyhedron $X$ is a set $Z\subset X$
such that $Z$ meets every maximal simplex $\Delta$, relative to a
triangulation $(K,\theta)$ (hence any,) in set whose pre-image
under $\theta$ has $n-$dimensional Lebesgue measure $0$,
$n=dim\Delta$. Note that 'almost everywhere' (a.e.) means
everywhere exept in some null set.

A {\it Riemannian polyhedron} $X=(X,g)$  is defined as a Lip
polyhedron $X$ with a specified triangulation $(K,\theta)$ such
that K is simplicial complex endowed with a covariant bounded
measurable Riemannian metric tensor $g$, satisfying the
ellipticity condition below. In fact, suppose that $X$ has
homogeneous dimension $n$ and choose a measurable Riemannian
metric $g_\Delta$ on the open euclidean $n-$simplex
$\theta^{-1}(\Delta^o)$ of $K$. In terms of euclidean coordinates
$\{x_1,...,x_n\}$ of points $x=\theta^{-1}(p)$, $g_\Delta$ thus
assigns to almost every point $p\in \Delta^o$ (or $x$), an
$n\times n$ symmetric positive definite matrix $g_\Delta =
(g_{ij}^\Delta(x))_{i,j=1,...,n}$ with measurable real entries and
there is a constant $\Lambda_\Delta
>0$ such that (ellipticity condition):
$$\Lambda_\Delta^{-2}\sum_{i=0}^{i=n}(\xi^i)^2\le \sum_{i,j}
g^\Delta_{ij}(x)\xi^i\xi^j\le\Lambda_\Delta^2\sum_{i=0}^{i=n}(\xi^i)^2$$
for $a.e.$ $x\in\theta^{-1}(\Delta^o)$ and every
$\xi=(\xi^1,...,\xi^n) \in \Bbb R^n$. This condition amounts to
the components of $g_\Delta$ being bounded and it is independent
not only of the choice of the euclidean frame on
$\theta^{-1}(\Delta^o)$ but also of the chosen triangulation.

For simplicity of statements we shall sometimes require that,
relative to a fixed triangulation $(K,\theta)$ of Riemannian
polyhedron $X$ (uniform ellipticity condition),
$$\text{$\Lambda$ $:=$
sup$\{\Lambda_\Delta:\Delta$ is simplex of $X\}<\infty$ .}$$

A Riemannian polyhedron $X$ is said to be admissible if for a
fixed triangulation $(K,\theta)$ (hence any) the Riemannian
simplicial complex $K$ is admissible.

We underline that (for simplicity) the given definition of a
Riemannian polyhedron $(X,g)$ contains already the fact (because
of the definition above of the Riemannian admissible complex) that
the metric $g$ is {\it continuous} relative to some (hence any)
triangulation (i.e. for every maximal simplex $\Delta$ the metric
$g_\Delta$ is continuous up to the boundary). This fact is some
times in the literature taken off. The polyhedron is said to be
simpexwise smooth if relative to some triangulation $(K,\theta)$
(and hence any), the complex $K$ is simplexwise smooth. Both
continuity and simplexwise smoothness are preserved under
subdivision.

In the case of a general bounded measurable Riemannian metric $g$
on $X$, we often consider, in addition to $g$, the {\it euclidean
Riemannian metric} $g^e$ on the Lip polyhedron $X$ with a
specified triangulation $(K,\theta)$. For each simplex $\Delta$,
$g^e_\Delta$ is defined in terms of euclidean frame on
$\theta^{-1}(\Delta^o)$ as above by unitmatrix $(\delta_{ij})$.
Thus $g^e$ is by no means covariantly defined and should be
regarded as a mere reference metric on the triangulated polyhedron
$X$.

Relative to a given triangulation $(K,\theta)$ of an
$n-$dimensional Riemannian polyhedron $(X,g)$ (not necessarily
admissible), we have on $X$ the distance function $e$ induced by
the euclidean distance on the euclidean space $V$ in which $K$ is
affinely Lip embedded. This distance $e$ is not intrinsic but it
will play an auxiliary role in defining an equivalent distance
$d_X$ as follows:

Let $\frak Z$ denote the collection of all null sets of $X$. For
given triangulation $(K,\theta)$ consider the set $Z_K\subset
\frak Z $ obtained from $X$ by removing from each maximal simplex
$\Delta$ in $X$ those points of $\Delta^o$ which are Lebesgue
points for $g_\Delta$. For $x,y \in X$ and any $Z\in \frak Z$ such
that $Z\subset Z_K$ we set:
$$\text{$d_X(x,y)=\sup \Sb {Z\in \frak Z}\\{Z\supset Z_K}\endSb
 \inf \Sb {\gamma}\\{\gamma(a)=x, \gamma(b)=y}\endSb \{ L_K(\gamma)$: $\gamma$ is Lip
continuous path and transversal to $Z\}$,}$$ where $L_K(\gamma)$
is de the length of the path $\gamma$ defined as:
$$\text{$L_K(\gamma)= \sum \Sb {\Delta\subset X}\\{}\endSb
 \int_{\gamma^{-1}(\Delta^o)} \sqrt{(g_{ij}^\Delta \circ \theta^{-1} \circ \gamma)
 \dot{\gamma}^i \dot{\gamma}^j } $, the sum is over all simplexes meeting $\gamma$.}$$

 It is shown in [12] that the distance $d_X$ is
 intrinsic, in particular it is independent of the chosen triangulation
 and it is equivalent to the euclidean distance $e$ (due to the Lip
 affinely and homeomorphically embedding of $X$ in some euclidean space $V$).

\enddefinition

\subhead { 4. Energy of maps}\endsubhead

The concept of energy in the case of a map of Riemannian domain
into an arbitrary metric space $Y$ was defined and investigated by
Korevaar and Shoen [21]. Later this concept was extended by Eells
and Fuglede [12] to the case of map from an admissible Riemannian
polyhedron $X$ with simplexwise smooth Riemannian metric. Thus,
The energy $E(\varphi)$ of a map $\varphi$ from $X$ to the space
$Y$ is defined as the limit of suitable approximate energy
expressed in terms of the distance function $d_Y$ of $Y$.

 It is shown in [12] that the maps $\varphi : X\flech Y$ of
 finite energy are precisely those quasicontinuous (i.e.
 has a continuous restriction to closed sets, whose complements
 have arbitrarily small capacity, (cf. [12] page 153) whose restriction
 to each top dimensional simplex of $X$
 has finite energy in the sense of Korevaar-Schoen, and
 $E(\varphi)$ is the sum of the energies of these restrictions.

 Just now, let $(X,g)$ be an admissible $m-$dimentional
 Riemannian polyhedron with simplexwise smooth Riemannian metric.
 It is not required that $g$ is continuous across lower
 dimensional simplexes. The target $(Y,d_Y)$ is an arbitrary
 metric space.

 Denote $L^2_{loc}(X,Y)$ the space of all $\mu_g-$mesurable
 ( $\mu_g$ the volume measure of $g$)
 maps $\varphi :X\flech Y$ having separable essential range and
 for which the map $d_Y(\varphi (.),q)\in L^2_{loc}(X,\mu_g)$
 (i.e. locally $\mu_g-$squared integrable)for some point $q$
 (hence by triangle inequality for any point). For $\varphi,\psi \in
 L^2_{loc}(X,Y)$ define their distance $D(\varphi,\psi)$ by:
 $$D^2(\varphi,\psi)= \int_X d_Y^2(\varphi (x),\psi(x)) d\mu_g(x).$$
 Two maps $\varphi,\psi \in L^2_{loc}(X,Y)$ are said to be {\it
 equivalent} if $D(\varphi,\psi)=0$, i.e. $\varphi(x)=\psi(x)$
 $\mu_g-$a.e. If the space $X$ is compact then $D(\varphi,\psi)<\infty$
 and $D$ is a metric on $L^2_{loc}(X,Y)=L^2(X,Y)$ and complete if
 the space $Y$ is complete [21].

 The {\it approximate energy density} of the map $\varphi\in
 L^2_{loc}(X,Y)$ is defined for $\epsilon >0$ by:
 $$e_\epsilon(\varphi)(x)=\int_{B_X(x,\epsilon)}\frac{d_Y^2(\varphi(x),\varphi(x'))}
 {\epsilon^{m+2}}d\mu_g(x').$$
 The function $e_\epsilon(\varphi)\ge 0$ is locally
 $\mu_g-$integrable.

 The {\it energy} $E(\varphi)$ of a map $\varphi$ of class
 $L^2_{loc}(X,Y)$ is:
 $$E(\varphi)=\sup_{f\in
 C_c(X,[0,1])}(\limsup_{\epsilon\rightarrow 0}\int_X f
 e_\epsilon(\varphi) d\mu_g),$$ where $C_c(X,[0,1])$ denotes the
 space of continuous functions from $X$ to the interval $[0,1]$
 with compact support.

A map $\varphi: X\flech Y$ is said {\it locally of finite energy},
and we write $\varphi \in W^{1,2}_{loc}(X,Y)$, if
$E(\varphi|U)<\infty$ for every relatively compact domain
$U\subset X$, or equivalently if $X$ can be covered by domains
$U\subset X$ such that $E(\varphi|U)<\infty$.

For example (cf. [12] lemma 4.4), every Lip continuous map
$\varphi : X \flech Y$ is of class $W^{1,2}_{loc}(X,Y)$. In the
case when $X$ is compact $W^{1,2}_{loc}(X,Y)$ is denoted
$W^{1,2}(X,Y)$ the space of all maps of finite energy.

$W^{1,2}_c(X,Y)$ denotes the linear subspace of $W^{1,2}(X,Y)$
consisting of all maps of finite energy of compact support in $X$.

We denote the closure of the space $Lip_c(X)$ (the space of
Lipschits continuous functions with compact supports) in the space
$W^{1,2}(X)$, $W_0^{1,2}(X)$.

We can show  (cf. [12] theorem 9.1) that a map $\varphi \in
L^2_{loc}(X)$ is locally of finite energy iff there is a function
$e(\varphi)\in L^1_{loc}(X)$, named {\it energy density} of
$\varphi$, such that (weak convergence):
$$\text{ $\lim_{\epsilon\rightarrow 0}\int_X f e_\epsilon (\varphi) d\mu_g
=\int_X f e(\varphi) d\mu_g$,  for each $f\in C_c(X)$.}$$


                                       \enddocument
\end